\documentclass[centertags,leqno]{article}

\usepackage{amsmath}
\usepackage{amssymb}
\usepackage{latexsym}
\usepackage{amsxtra}
\usepackage{amscd}
\usepackage{theorem}
\usepackage{xypic}

\theoremstyle{plain}

{\theorembodyfont{\slshape}

        \newtheorem{thm}{Theorem}[section]
        \newtheorem{cor}[thm]{Corollary}
        \newtheorem{lem}[thm]{Lemma}
        \newtheorem{prop}[thm]{Proposition}

}

{\theorembodyfont{\rmfamily}

        \newtheorem{defn}[thm]{Definition}
        
        \newtheorem{rem}[thm]{Remark}

        \newtheorem{exa}[thm]{Example}

}

\newcommand{\proof}{{\bf Proof:\ }}

\newcommand{\Endproof}{\hspace*{\fill} $\Box$ \vspace{1ex} \noindent }

\makeatletter
\renewcommand{\subsection}{\@startsection{subsection}{2}%
        {\z@}{-3.25ex plus -1ex minus-.2ex}{-1em}{\bf}}
\makeatother

\newcommand{\PP}{\mathbb{P}}

\newcommand{\CC}{\mathbb{C}}

\newcommand{\HH}{\mathcal{H}}

\newcommand{\ord}{{\rm ord}}

\title{Existence of covers with fixed ramification in positive characteristic}
\author{Irene I. Bouw and Leonardo Zapponi}
\date{}
\begin{document}
\maketitle

\begin{abstract}
We discuss two elementary constructions for
covers with fixed ramification in positive characteristic. As an application, we compute
the number of certain classes of covers between projective lines
branched at $4$ points and obtain information on the structure of the
Hurwitz curve parametrizing these covers.
\end{abstract}

\section{Introduction}
In this note, we consider the question of determining the number of
covers between projective lines in positive characteristic with
specified ramification data and fixed branch points. The ramification
data considered are the degree of the cover, together with a list of
the ramification indices in the fibers of the branch points. Over an
algebraically closed field of characteristic zero, it is in principal
 possible to solve this problem by Riemann's Existence Theorem. Namely, the number of covers can be expressed as the cardinality of a finite set,
which can be explicitly constructed in concrete cases. In particular, this approach
shows that the number of covers is finite and does not depend on the
position of the branch points. 

In positive characteristic, the situation is drastically different. For example, the
number of covers with fixed ramification  depends on the position
of the branch points. Moreover, if the characteristic $p$ divides one
of the ramification indices, the number of covers is in general
infinite. There are only few general results on the number of covers
in this situation (we refer to \cite{BO} for an overview).

 The work of Osserman (\cite{Osserman1}, \cite{Osserman2}, \cite{LO})
 suggests that a particularly nice case to look at is that of covers
 $f:\PP^1\to \PP^1$ of degree $d$ which are ramified at $r$ points
 $x_1, \ldots, x_r$ with $f(x_i)$ pairwise distinct (the
 so-called {\sl single-cycle} case).  We write $h(d; e_1, e_2, e_3,
 \ldots, e_r)$ for the number of single-cycle covers with fixed branch
 locus over $\CC$, where $e_i$ is the ramification index of
 $x_i$; this number is called the {\sl Hurwitz number}.

Let $k$ be an algebraically closed field of positive characteristic $p$.  We
only consider covers $f:\PP^1_k\to \PP^1_k$ in the tame and single-cycle
case. We denote by $h_p(d; e_1, e_2, e_3, \ldots,
e_r)$ the maximal number of covers with fixed branch locus, where
the maximum is taken over all possible branch loci. This number is called the
$p$-{\sl Hurwitz number}. Since $p\nmid e_i$ for all $i$, this number
is finite and does not depend on $k$.  It can be shown that there the maximum is attained if the branch locus belongs to a dense open subset $U\subset (\PP^1_k)^r\setminus \Delta$. Here $\Delta$
is the fat diagonal.

We start by summarizing the results on the number of covers with fixed
branch locus in the single-cycle case for $r\in \{3,
4\}$. In \cite{LO}, F. Liu and B. Osserman give a closed formula for the number
of such covers in characteristic zero. In \cite{Osserman1} and 
\cite{Osserman2}, B. Osserman determines the $p$-Hurwitz number $h_p(d;
e_1, e_2, e_3)$ using linear series. In \cite{BO} the number $h_p(p;
e_1, e_2, e_3, e_4)$ is computed. This last case is substantially more
difficult and he proof relies on the theory of stable reduction of covers. 

In this note, we also consider covers $f:\PP^1_k\to \PP^1_k$ of
ramification type $(d; e_1, e_2, e_3, e_4)$. In contrast with the
situation in \cite{BO}, the degree $d$ is not fixed.  We consider two
elementary constructions, which yield previously unknown results on
some $p$-Hurwitz numbers $h_p(d; e_1, e_2, e_3, e_4)$. Both
constructions were known before and can be found for example in
\cite{Osserman1}. However, the implications for the $p$-Hurwitz
numbers have not been fully exploited. As an additional result, we
obtain rather complete information on the structure of the Hurwitz
curve, parameterizing covers of the type considered, in positive
characteristic. These are the first such results.

The first result deals with the case $1< e_i< p$ and $e_4=p-1$. In
this situation, we compute the $p$-Hurwitz number $h_p(d; e_1, e_2, e_3, e_4)$. We can even obtain
something stronger, namely an explicit description of the Hurwitz curve
$\HH_p(d; e_1, e_2, e_3, e_4)$ parameterizing all covers of type $(d;
e_1, e_2, e_3, e_4)$. This yields, in particular, not only a formula for
all covers with generic branch locus, but also exactly describes the values for which  the number of covers drops. As far as we know, this is the first nontrivial example of a complete description of Hurwitz curves in positive characteristic. Other papers on Hurwitz curves in positive characteristic (eg
\cite{crelle}, \cite{meta} and \cite{BO}) do not yield such description.
We refer to \S\ \ref{multconstsec} for the precise
statement of the result.

The second result considers the case  $e_1>p$ and $2<e_2, e_3,
e_4<p$. Here, we relate the $p$-Hurwitz numbers $h_p(d; e_1,
e_2, e_3, e_4)$ and $h_p(d; e_1, e_2, e_3$-$e_4)$. In general none of these
two numbers is known. However, we illustrate in a concrete example the kind
of  results obtained  using this method.

\section{Notation and basic results}
Let $k$ be an algebraically closed field. We consider tamely ramified
covers $f:X\to \PP^1_k$ between smooth projective curves. Let
${\boldsymbol x}=(x_1=0, x_2=1, x_3=\infty, x_4, \ldots, x_r)$
be the branch points of $f$, which we consider to be ordered. We firstly associate to $f$ a ramification invariant.

\begin{defn}\label{raminvdef} Let $f:X\to \PP^1_k$ be a tamely ramified cover
 as above. Denote by $d$ the degree of $f$. For every $i$, let $(e_{i,
   1}, \ldots e_{i, n_i})$ be the partition of $d$ corresponding to
 the ramification indices of the points in the fiber
 $f^{-1}(x_i)$. This partition defines a conjugacy class $C_i$ of the
 symmetric group $S_d$. The {\sl ramification type} of $f$ is defined
 as the datum ${\boldsymbol C}=(d; C_1, \ldots, C_r)$. If $C_i$ is the class of a single cycle (i.e.\ exactly one $e_{i,j}$ is different from $1$) we simply write $C_i=e_i$, where $e_i$
 is the length of the cycle. If $p$ is a prime number, we say that ${\boldsymbol C}$ is $p$-{\sl tame} if all $e_{i, j}$ are relatively prime to $p$.
\end{defn}

Since we assume that $f$ is tamely ramified, the Riemann--Hurwitz
formula states that
\begin{equation}\label{RHeq}
\sum_{i, j} e_{i, j}-\sum_i n_i=2g(X)-2+2d.
\end{equation}
Since the genus $g(X)$ of $X$ only depends on the ramification type of the cover, we
sometimes denote it by $g({\boldsymbol C})$ and refer to it as the
{\sl genus of the ramification type}.  A datum $(d; C_1, \ldots,
C_r)$, where $C_i=(e_{i, 1}, \ldots e_{i, n_i})$ is a partition of $d$
and $g({\boldsymbol C})$ an integer, is called a {\sl ramification
  type}. A ramification type such that $g({\boldsymbol C})=0$ is
called a {\sl genus}-$0$  type.

  There exists
several variants of these definitions. For example, in some cases it
makes sense to include the Galois group of the Galois closure of $f$
into the definition. In this note, we mostly consider the case that
each $C_i=e_i$ is the conjugacy class of a single cycle. In this case,
the Galois group is typically $S_d$ or $A_d$, with few exceptions in
small degree.

Two covers
 $f_i:X_i\to \PP^1$  are considered
 {\sl isomorphic} if there exists an isomorphism $\iota:X_1\to X_2$ making
\[
\xymatrix{
X_1\ar[rr]^{\sim}_{\iota}\ar[dr] &&X_2\ar[dl]\\
&\PP^1&
}
\]
commutative.
In particular, both covers have the same branch locus and the same ramification type.

 The number of isomorphism classes of
 covers with a given ramification type $(d; C_1, \ldots, C_r)$ and fixed
 branch points ${\boldsymbol x}$ is finite, since we only consider tame ramification.  We first
 consider the classical characteristic zero case, for which the number of
 isomorphism classes of covers does not depend on the position of the branch
 points. It follows from  Riemann's Existence
 Theorem that this number,  called the {\sl Hurwitz number} and denoted by
 $h(d; C_1, \ldots, C_r)$, is the cardinality of the set
\[
\left\{(g_1, \ldots, g_r)\in C_1\times \cdots C_r\mid \langle
g_i\rangle\subset S_d \text{ transitive},\, \prod_i g_i=e\right\}/\sim,
\]
where $\sim$ denotes uniform conjugacy by the group $S_d$. The
condition that the $g_i$ generate a transitive subgroup of $S_d$
guarantees that the corresponding cover $f$ is connected.

 Let ${\boldsymbol C}:=(d; C_1, \ldots, C_r)$ be a ramification
 type. We denote by $\HH_k({\boldsymbol C})$ the Hurwitz space
 parameterizing isomorphism classes of covers $f:Y\to \PP^1_k$ with
 ramification type ${\boldsymbol C}$ defined over $k$.

Let $\pi: \HH_k({\boldsymbol C}) \to (\PP^1_k)^{r-3}\setminus \Delta$
be the natural map $[f]\mapsto {\boldsymbol x}$ which sends the class
of a cover to its branch locus. Here, $\Delta:=\{{\boldsymbol x}\mid
  x_i= x_j \text{ for some }i\neq j\}$ is the fat diagonal.  In
  characteristic zero, this map is finite and flat. Its degree is
  exactly the Hurwitz number $h(d; C_1, \ldots, C_r)$. Moreover, the
  map $\pi$ is unramified. We remark that the Hurwitz space may not be be connected.

Now assume that the characteristic $p$ of $k$ is positive. Then the number of
covers of given ramification type may depend on the position of the
branch points. The {\sl $p$-Hurwitz number} $h_p({\boldsymbol C})$ is
defined as the number of isomorphism classes of covers of ramification
type ${\boldsymbol C}$ for which the branch locus ${\boldsymbol x}$ is
generic, in the sense that it corresponds to the generic point of
$(\PP^1)^{r-3}\setminus \Delta$.  The $p$-Hurwitz number is also the
maximum number of covers of given type as the branch locus $\boldsymbol x$ varies.

The following well-known lemma gives some information on the
Hurwitz number in this context. Part (a) follows from the fact that every
tame cover in characteristic $p$ lifts to characteristic zero. Part
(b) is a consequence of the isomorphism between the prime-to-$p$ part of the fundamental group
$\pi^{(p)}(\PP^1_k\setminus{\boldsymbol x}, \ast)$ and the
prime-to-$p$ part of the fundamental group of the complement of $r$
points on $\PP^1$ in characteristic zero (\cite{SGA1}).

\begin{lem}\label{p-hurwitzlem}
Let $k$ be an algebraically closed field of characteristic $p>0$.
\begin{itemize}
\item[(a)] The $p$-Hurwitz number $h_p({\boldsymbol C})$ only depends
  on $p$, and not on the field $k$.
\item[(b)] We have $h_p({\boldsymbol C})\leq h({\boldsymbol C})$ with
  equality if $d<p$.
\end{itemize}
\end{lem}

Note that the difference $h({\boldsymbol C})-h_p({\boldsymbol C})$
corresponds to the number of covers in characteristic zero which have
generic branch locus and bad reduction to characteristic $p$. We call
this number sometimes the {\sl bad degree} of the ramification type
${\boldsymbol C}$ if $p$ is clear from the context.

Suppose that ${\boldsymbol C}=(d; C_1, \ldots, C_r)$ is a genus-$0$
ramification type. In this case the degree $d$ of a cover of type
${\boldsymbol C}$ is determined by the conjugacy classes $C_i$ via the
Riemann--Hurwitz formula (\ref{RHeq}). For convenience, we may therefore drop $d$ from
the notation. We write $h^0( C_1,
\ldots, C_r)=h(d; C_1, \ldots, C_r)$ and $h_p^0( C_1, \ldots,
C_r)=h_p(d; C_1, \ldots, C_r)$ for the Hurwitz numbers in the
genus-$0$ case. 

The following lemma describes the Hurwitz numbers in the $3$-point
single-cycle case under the assumption that the genus of the
ramification type is zero.

\begin{lem}\label{3ptlem}
Let ${\boldsymbol C}=(d; e_1, e_2, e_3)$ be a ramification type with
$g({\boldsymbol C})=0$. Then
\begin{itemize}
\item[(a)] $h^0( e_1, e_2, e_3)=1$.
\item[(b)] Assume additionally that $e_1, e_2<p$. Then $h_p^0(e_1,
  e_2, e_3)=1$ if and only if $d<p$ and $0$ otherwise.
\end{itemize}
\end{lem}

\proof Part (a) is an elementary calculation using the combinatorial
description of the Hurwitz number given above (see for example
\cite[Lemma 2.1]{LO}). Part (b) is proved by B. Osserman (\cite[Cor
  2.5]{Osserman2}) using linear series.  \Endproof

B. Osserman proves a stronger version of Lemma \ref{3ptlem}, calculating
$h_p(d; e_1, e_2, e_3)$ in the situation of Lemma \ref{3ptlem}, but
without the assumption that $e_1, e_2<p$. We do not recall the full
result here, since it its formulation if quite involved. The following
proposition, proved in \cite{LO}, is a generalization of Lemma
\ref{3ptlem}.(a) to the case of $4$ branch points.

\begin{prop}\label{4ptprop}
Let ${\boldsymbol C}=(d; e_1, e_2, e_3, e_4)$ be a ramification type with
$g({\boldsymbol C})=0$. Then
\begin{itemize}
\item[(a)] $h^0({\boldsymbol C})=\min_i(e_i(d+1-e_i))$.
\item[(b)] The Hurwitz curve $\HH_{\CC}({\boldsymbol C})$ is
  connected.
\end{itemize}
\end{prop}

\begin{rem} Assume that ${\boldsymbol C}=(d;e_1, \ldots, e_r)$ is a single-cycle
 ramification type with genus $g({\boldsymbol C})=0$.  Write
 ${\boldsymbol x}=(x_1=0, x_2=1, x_3=\infty, x_4 \ldots, x_r)$
 for the branch points and ${\boldsymbol y}=(y_1, \ldots y_r)$ for the
 ramification points, where $f(y_i)=x_i$. Up to isomorphism, the
 associated cover $f:Y\simeq \PP^1_k\to X\simeq\PP^1_k$ may be
 normalized such that $y_1=0, y_2=1$ and $y_3=\infty$. If this is the
 case, we say that $f$ is {\sl normalized}. Note that any isomorphism
 class contains a unique normalized representative. Assuming $f$ is
 normalized, we may therefore regard $f$ as element of $k(T)$, where
 $x$ is a coordinate on $Y\simeq \PP^1_k$ with $T(i)=i$ for $i\in \{0,
 1, \infty\}$.
\end{rem}

\section{An elementary construction}\label{multconstsec}
We start by recalling an elementary construction due to B. Osserman
\cite[Lemma 5.2]{Osserman1}. Let $1<e_1, \ldots, e_r<p$. Osserman's
result establishes a bijection between maps $f:\PP^1_k\to \PP^1_k$
with ramification indices $(e_i)$ and maps $h:\PP^1_k\to \PP^1_k$ with
ramification indices $p-e_1, e_2, \ldots, e_{r-1}, p-e_r$. The maps
$f$ and $h$ need not have the same degree. Although B. Osserman does not
explicitly consider this, the construction also works when either
$e_1$ or $e_r$ equals $p-1$ in which case the number of ramification
points of $f$  differs from that of $h$.

In the rest of this section, we assume that $r=4$. The following lemma
is a normalized version of the result of Osserman. In our
version, we make sure that the ramification points map to distinct
points. This allows to deduce a statement on Hurwitz numbers in
characteristic $p>0$.

Let $k$ be an algebraically closed field of characteristic $p>0$. We
fix a genus-$0$ ramification type ${\boldsymbol C}=(d; e_1, e_2, e_3,
p-1)$ and a branch locus ${\boldsymbol x}=(x_1=0, x_2=1, x_3=\infty,
x_4=:\lambda)$. We assume that $1<e_i<p$ for all $i$. Note that the degree of a cover of
type ${\boldsymbol C}$ is equal to $d=(e_1+e_2+e_3+p-3)/2$.

\begin{lem}\label{multconstlem}
\begin{itemize}
\item[(a)] Let $f:\PP^1_k\to \PP^1_k$ be a normalized cover of type
  ${\boldsymbol C}=(d; e_1, e_2, e_3, p-1)$ and branch locus
  ${\boldsymbol x}=(x_1=0, x_2=1, x_3=\infty, x_4=:\lambda)$. We
  denote the unique ramification point of $f$ above $\lambda$ by $\mu$.
  Then the cover $h:\PP^1_k\to \PP^1_k$ defined by
\[
g(y)=\frac{(y-\mu)^p}{f-\lambda}, \quad h(y):=\frac{g-g(0)}{g(1)-g(0)}
\]
is a normalized cover of type $\tilde{{\boldsymbol
  C}}:=(\tilde{d};e_1, e_2, p-e_3)$ and branch locus
$ \tilde{{\boldsymbol x}}=(x_1=0, x_2=1,x_3=\infty)$. Here
$\tilde{d}=d+1-e_3=(e_1+e_2-e_3+p-1)/2$.
\item[(b)] Conversely, choose an element $\mu\in
  \PP^1_k\setminus\{0,1, \infty\}$ with $h(\mu)\neq 0, 1, \infty,
  \mu^p$. Suppose given a normalized cover $h:\PP^1_k\to \PP^1_k$ of
  ramification type $\tilde{{\boldsymbol C}}:=(\tilde{d};e_1, e_2,
  p-e_3)$.  Then the cover $f:\PP^1_k\to
  \PP^1_k$ defined by
\[
w(y)=\frac{(y-\mu)^p}{h-h(\mu)}, \quad f(y)=\frac{w-w(0)}{w(1)-w(0)}
\]
is a normalized cover of type ${\boldsymbol C}=(d; e_1, e_2, e_3,
p-1)$ with branch locus ${\boldsymbol x}=(0, 1, \infty,
\lambda:=\mu^p(1-h(\mu)/(\mu^p-h(\mu))$.
 \item[(c)] The constructions of (a) and (b) are inverse to each other.
\end{itemize}
\end{lem}

\proof Let $f$ be as in  (a). The statement on the
ramification indices of $h$ follows immediately from the definition of
$h$. To prove the statement of the ramification type, it remains to
show that $h(0), h(1), h(\infty)$ are pairwise distinct.  The condition
$h(0)=h(1)$ is equivalent to $\mu^p=\lambda$. If this were the case,
the cover $h$ would only have $2$ branch points, with ramification
type ${\boldsymbol C}^\ast=(\tilde{d}; e_1$-$e_2$-$1\cdots$-$1,
(p-e_3)$-$1\cdots$-$1)$. It is well-known that such a cover does not
exist. This shows that $h$ has the claimed ramification type. The
statement on the degree follows from the observation that $h$
 is a genus-$0$ cover.

Let $h$ be as in (b), and write $h(y)=y^{e_1}h_1/h_2$, where
$\deg(h_1)=\tilde{d}-e_1$ and $\deg(h_2)=\tilde{d}-p+e_3$. Note that
the definition of $w$ may be rewritten as
\[
w(y)=\frac{(y-\mu)^{p-1}h_2}{h_3},
\]
where $h_3:=(y^{e_1}h_1-h(\mu)h_2)/(y-\mu)$ is a polynomial of
degree $\tilde{d}-1$. It follows that $f$ maps the ramification points
$0, 1, \infty$ to the points $0, 1, \infty$,  respectively.

Define $\lambda:=f(\mu)$. An easy computation leads to the identity
\[
\lambda=f(\mu)=\frac{-w(0)}{w(1)-w(0)}=\frac{\mu^p(1-h(\mu))}{\mu^p-h(\mu)}.
\] 
The choice of $\mu$  implies that $f(0), f(1), f(\infty), f(\mu)$
 are pairwise distinct.  
Therefore (b) follows similarly to (a).

Part (c) is an easy verification. We prove one direction and leave the
other as an exercise. Let $h:\PP^1_k\to \PP^1_k$ be a normalized cover
of type $\tilde{{\boldsymbol C}}$, and define $f$ as in the statement
of (b). The definition of $\lambda$ implies  that
\[
f-\lambda=\frac{w}{c}=\frac{(y-\mu)^p}{c(h-h(\mu))},
\]
where have set $c:=w(1)-w(0)$. The cover $g$ associated with
 $f$ in  (a) satisfies $g=c(h-h(\mu))$. Since $h$ is
 normalized, we conclude that the unique normalized
polynomial associated with $g$ is again $h$.  \Endproof

%Write $\HH_p({\boldsymbol
%  x}; {\boldsymbol C})$ for the set of isomorphism classes of
%normalized covers of type ${\boldsymbol C}$ with branch locus
%${\boldsymbol x}$. We may identify this set with the fiber
%$\pi^{-1}({\boldsymbol x})\subset \HH_p({\boldsymbol C})$.

\begin{prop}\label{hurwitznrprop} Let $p$ be an odd prime number and 
$1<e_1, e_2, e_3<p$ integers such that $E:=e_1+e_2+e_3$ is even and
  $e_3\neq p-1$. Put $d=(E+p-3)/2$, i.e.\ ${\boldsymbol C}:=(d; e_1,
  e_2, e_3, p-1)$ is a genus-$0$ ramification type.
\begin{itemize}
\item[(a)] The Hurwitz number $h_p^0({\boldsymbol C})$ is positive if
  and only if
\begin{equation}\label{phposeq}
p+1\leq E\leq p-1+2\min_i(e_i).
\end{equation}
\item[(b)] If the equivalent conditions of (a) are satisfied, we have
  that
\[
h_p^0({\boldsymbol C})=\frac{1}{2}(3p-1-E).
\]
\end{itemize}
\end{prop}

\proof Write $\tilde{\boldsymbol{C}}=(\tilde{d}; e_1, e_2, p\text{-}
e_3)$, where $\tilde{d}=d=1-e_3$.  Lemma \ref{multconstlem} implies
that $h_p^0({\boldsymbol C})$ is positive if and only if
$h_p^0(\tilde{\boldsymbol{C}})$ is positive.  A necessary and
sufficient condition for $\tilde{\boldsymbol C}$ to be a ramification
type with three branch points is that
\begin{equation}\label{ramtypeeq}
1<e_1, e_2, p-e_3\leq \tilde{d}.
\end{equation}
In our situation, the conditions (\ref{ramtypeeq}) 
can be rewritten as
\[
\begin{cases}
e_1, e_2&>1,\\
e_3&<p-1,\\
p+1&\leq e_1+e_2+e_3,\\
e_1+e_2+e_3&\leq p-1+2e_1,\\
e_1+e_2+e_3&\leq p-1+2e_2.
\end{cases}
\]

Assume that these conditions are satisfied. Lemma \ref{3ptlem}.(b) implies that
 $h_p^0(e_1, e_2, p-e_3)$ is positive if and only if
\begin{equation}\label{degeq}
\tilde{d}=(e_1+e_2-e_3+p-1)/2<p.
\end{equation}
Combining these inequalities immediately yields (a).

To prove (b), we assume that $h_p^0(e_1, e_2, p-e_3)$ is positive,
i.e.\ that (\ref{phposeq}) holds. Lemma \ref{p-hurwitzlem}.(b) implies
that $h^0(e_1, e_2, p-e_3)$ is positive as well, hence $h^0(e_1, e_2,
p-e_3)=1$ (cf. Lemma \ref{3ptlem}.(a)). Applying Lemma
\ref{p-hurwitzlem}.(b) and using the fact that $\tilde{d}<p$, we 
therefore obtain the identity $h_p^0(e_1, e_2, p-e_3)=1$.

Denote by $h:\PP^1_k\to \PP^1_k$ the unique normalized cover of type
$(e_1, e_2, p-e_3)$. Choose $\mu$ as in Lemma \ref{multconstlem}.(b),
and define $f$ and $\lambda$ as the statement of Lemma
\ref{multconstlem}.(b). Consider the birational map 
\begin{equation}\label{lambdaeq}
\mu\mapsto \lambda(\mu)=\frac{\mu^p(1-h(\mu))}{\mu^p-h(\mu)}, \quad
\PP^1_k\to \PP^1_k.
\end{equation}
Obviously, the degree $\deg(\lambda)$ of this map equals the
$p$-Hurwitz number $h_p^0(e_1, e_2, p-e_3)$. We compute
$\deg(\lambda)$ by computing the divisor of this map.

Since $h$ is normalized, it follows that
\[
\ord_{\mu=0}(\lambda)=p-e_1, \quad \ord_{\mu=1}(\lambda)=e_2-e_2=0,
\quad \ord_{\mu=\infty}(\lambda)=-(p-e_3)<0.
\]
Moreover, $\lambda$ has $\tilde{d}-e_2$ simple zeros, which are
different from $\mu=0,1,\infty$. Note that $\mu^p-h(\mu)$ has
$\tilde{d}-p+e_3$ poles different from $\mu=\infty$, which all have
multiplicity one. However, these are exactly the (simple) poles of
$1-h(\mu)$ which are different from $1$, hence these do not yield
zeros of $\lambda$. We conclude that
\[
\deg(\lambda)=p-e_1+\tilde{d}-e_2=(3p-1-(e_1+e_2+e_3))/2,
\]
which proves (b).
\Endproof

\begin{rem} B. Osserman gives a similar 
statement  (\cite[Cor.\ 8.1]{Osserman1}), counting covers with $4$ ramification points. His result is
more general, since it does not require one of the
ramification indices to be $p-1$.  However, he fixes the ramification
points of the cover, rather than the branch points. Therefore his
result does not compute Hurwitz numbers in characteristic $p>0$.

Computing $p$-Hurwitz numbers is in general more difficult than
counting covers with fixed ramification. Beside the classical result
from Lemma \ref{p-hurwitzlem}.(b), the only general result on
$p$-Hurwitz numbers is the main result of \cite{BO}, which
computes $h_p(p; e_1, e_2, e_3, e_4)$. That result relies on subtle
and deep results on the stable reduction of Galois covers. 
\end{rem}

The following corollary translates the statement of Proposition
\ref{hurwitznrprop} into a statement on the Hurwitz curve
$\HH_p({\boldsymbol C})$.

\begin{cor}\label{hurwitznrcor} Let ${\boldsymbol C}=(d; e_1, e_2, e_3, p-1)$
 be as in the 
statement of Proposition \ref{hurwitznrprop}. The Hurwitz curve
$\HH_p({\boldsymbol C})$ is connected.
\end{cor}

\proof The statement immediately follows from the proposition. Let
$\pi_p:\HH_p({\boldsymbol C})\to \PP^1_\lambda$ be the natural map which
sends a cover of type ${\boldsymbol C}$ to the branch point
$\lambda$. Then $\pi$ is birationally equivalent to the map
$\mu\mapsto \lambda$ described in the prof of that proposition.
\Endproof

Let ${\boldsymbol C}:=(d; e_1, e_2, e_3, p-1)$ be a
ramification type satisfying the equivalent conditions of Proposition
\ref{hurwitznrprop}.(a). Put $\tilde{{\boldsymbol C}}=(\tilde{d}; e_1,
e_2, p-e_3)$ and let $h:\PP^1_k\to \PP^1_k$ be the unique cover of type
$\tilde{{\boldsymbol C}}$ (compare with the proof of Proposition
\ref{hurwitznrprop}.(a)). We may write $h(y)=h_1(y)/h_2(y)$, where the
$h_i\in k[y]$ are relatively prime and satisfy the relations
$\deg(h_1)=\tilde{d}-e_1$ and
$\deg(h_2)=\tilde{d}-(p-e_3)=(e_1+e_2+e_3-p-1)/2$. 

It follows from Lemma \ref{p-hurwitzlem} that there exist finitely many values $\lambda\in \PP^1_\lambda\setminus\{ 0,1,\infty\}$ for
which the number of covers of ramification type ${\boldsymbol C}$ and
branch locus $(0, 1, \infty, \lambda)$ is strictly less than
$h_p({\boldsymbol C})$.  We let $\Sigma({\boldsymbol C})\subset
\PP^1_\lambda\setminus\{0,1,\infty\}$ be this exceptional set and
call it the {\sl supersingular locus}.

\begin{cor}\label{supersingcor}
With the above notation, we have
\[
\Sigma({\boldsymbol C})=\{y\in \PP^1_k\setminus\{0,1,\infty\}\mid h_2(y)=0\}.
\]
\end{cor}

\proof We recall that the construction of Lemma
\ref{multconstlem}.(b) works if and only if $h(\mu)\neq 0,1, \infty,
\mu^p$.  Moreover, equation (\ref{lambdaeq}) gives an expression of the fourth branch point $\lambda$ of $f$ as function of $\mu$. 

Assume that $h(\mu)=0$. Then (\ref{lambdaeq}) implies that either
$\mu=0$ or $h(\mu)=1$. By definition $0, 1\not\in \Sigma({\boldsymbol
  C})$. Therefore it suffices to consider the solutions of $h(\mu)=1$
with $\mu\neq 1$. We may write $h(\mu)-1=\mu^{e_2}\varphi(\mu)$, where
$\varphi(1)\neq 1$. Substituting this in (\ref{lambdaeq}) yields
\[
\lambda(\mu)=\frac{\mu^p\varphi(\mu)}{(\mu-1)^{p-e_2}-\varphi(\mu)}.
\]
In particular, it follows that $\lambda(\mu)=0$ if $\mu$ is a zero of
$\varphi$. Hence these zeroes are not contained in $\Sigma({\boldsymbol
  C})$. Similarly it follows that the solutions of $h(\mu)=1$ don't belong to $\Sigma({\boldsymbol C})$.

Assume that $h(\mu)=\infty$ and $\mu\neq \infty$, i.e.\ $h_2(\mu)=0$ according to the notation introduced above the statement of the corollary. We then have the identity $\lambda(\mu)=\mu^p$. Therefore $\mu\in
\Sigma({\boldsymbol C})$, since $\mu\neq 0,1,\infty$.

Finally, assume that $\mu=\mu^p$ and $\mu\not\in\{0,1,\infty\}$. Then
$\lambda=\infty$, hence this does not yields any new value.
\Endproof

\begin{exa} We illustrate the results of this section with two concrete 
examples.

(a) Let $p\geq 5$ be a prime, and consider the genus-$0$ ramification
type ${\boldsymbol C}=(d; 2,2,p-3, p-1)$. Note that the condition of
Proposition \ref{hurwitznrprop}.(a) is satisfied. Hence Proposition
\ref{hurwitznrprop}.(b) implies that $h_p({\boldsymbol C})=p-1$ and
Proposition \ref{4ptprop}.(a) leads to $h({\boldsymbol
  C})=\min(3(p-3), 2(p-2),p-1)=p-1$. We therefore find the equality $h({\boldsymbol
  C})=h_p({\boldsymbol C})$, so that all covers of this type with generic
branch locus have good reduction.

The unique normalized cover $h:\PP^1\to \PP^1$  of
type $(3; 2,2,3)$ is given by
\[
h(y)=3y^3-2y^2.
\]
Therefore 
\[
\lambda(\mu)=\frac{\mu^{p}(1+2\mu^3-3\mu^2)}{\mu^p+2\mu^3-3\mu^2}=
\frac{\mu^{p-2}(2\mu-1)}{\sum_{i=1}^{p-4}i
  \mu^{p-3-i}}.
\]
This confirms that the degree $\deg(\lambda)$ equals
$(3p-(e_1+e_2+e_3))/2=p-1$.

We have already remarked that all covers with generic branch locus
have good reduction to characteristic $p>0$. Arguing as in the proof
of \cite[Theorem 4.2]{meta}, one may deduce from this observation that
the map $\pi_p:\HH_p({\boldsymbol C})\to
\PP^1_\lambda\setminus\{0,1,\infty\}$ is finite.  However, in this
concrete example the finiteness of $\pi_p$ immediately follows from
Corollary \ref{supersingcor}.

(b) Next, we consider the case ${\boldsymbol C}=(p; 3,2,p-2, p-1)$
and ${\tilde{\boldsymbol C}}=(3; 3,2,2)$, again assuming that $p\geq
5$. The unique cover $h:\PP^1_k\to \PP^1_k$ of type
${\tilde{\boldsymbol C}}$ is given by $h(y)=y^3/(3y-2)$ and a direct computation leads to the expression
\[
\lambda(\mu)=\frac{\mu^{p-3}(-\mu^3+3\mu-2)}{3\mu^{p-2}-2\mu^{p-3}-1}.
\]
Dividing the numerator and the denominator by $(\mu-1)^2$, we find
$\deg(\lambda)=p-2$, which confirms Proposition
\ref{hurwitznrprop}.(b). As in Corollary \ref{supersingcor}, the supersingular values are the poles of $h$ different from $\infty$. In this concrete example, we find a unique value, namely $\mu=2/3$. The case $d=p$ has been considered in \cite{BO}. This result may also be deduced from
\cite[Remark 8.3]{BO} (note, however, that the proof relies on subtle
arguments involving stable reduction, which are only sketched in that
paper). As in the previous example, Proposition \ref{4ptprop} implies that $h({\boldsymbol
  C})=2(p-1)$ and Proposition \ref{hurwitznrprop}.(b) asserts that
$h_p({\boldsymbol C})=p-2$ We conclude that $h({\boldsymbol
  C})-h_p({\boldsymbol C})=p$, which confirms the main result of \cite{BO}.
\end{exa}

Fix a genus-$0$ ramification type ${\boldsymbol C}=(d; e_1, e_2, e_3,
p-1)$ which is $p$-tame, i.e.\ $p\nmid e_i$.  Recall that
$h({\boldsymbol C})-h_p({\boldsymbol C})$ denotes the ``bad degree''
of the ramification type. This is the number of covers with generic
branch locus which have bad reduction to characteristic $p$. 

\begin{prop}\label{baddegprop}
The notation being as above, assume that the minimum\\ $\min_{i\in
  \{1,2,3\}} e_i(d+1-e_i)$ is attained for $e_1$. This is not a
restriction, since\\ we may permute the branch points. Then, the bad degree
$h(\boldsymbol{C})-h_p(\boldsymbol{C})$ is\\ given by
\[
\begin{cases} 
0& \text{ if }d\leq p-1,\\
p(d+1-p)& \text{ if }p\leq d\leq p-2+e_1,\\
h({\boldsymbol C})=e_1(d+1-e_1)&\text{ otherwise}.
\end{cases}
\]
\end{prop}

\proof The first case immediately follows from Lemma
\ref{p-hurwitzlem}.  Assume that $p\leq d\leq p-2+e_1$. In this case, Propositions \ref{4ptprop}.(a) and
\ref{hurwitznrprop}.(a) assert that $h({\boldsymbol C})=(p-1)(d+2-p)$ and
$h_p({\boldsymbol C})\neq 0$. Statement (b) therefore follows from
Proposition \ref{hurwitznrprop}.(b).

For $d> p-2+e_1$, Proposition \ref{4ptprop}.(a) implies
that $h({\boldsymbol C})=e_1(d+1-e_1)$. Since $d>p-2+\min_i
e_i=p-2+e_1$ by assumption, we conclude from Proposition
\ref{hurwitznrprop}.(a) that $h_p({\boldsymbol C})=0$ and statement (c)
follows.  \Endproof

In the second case of Proposition \ref{baddegprop}, some covers have
good reduction while others have bad reduction. In the first
(resp.\ third) case all covers have good (resp.\ bad) reduction to
characteristic $p>0$.  The following corollary therefore follows
from Proposition \ref{baddegprop} and its proof. A similar phenomenon
occurs in the situation of \cite[Section 4]{meta}.

\begin{cor}\label{baddegcor}
Let ${\boldsymbol C}$ be as in Proposition \ref{baddegprop}, and
assume that $h({\boldsymbol C})\neq h_p({\boldsymbol C})\neq 0$. Then
the bad degree $h({\boldsymbol C})-h_p({\boldsymbol C})$ is divisible
by $p$.
\end{cor}

\section{A variant}\label{addconstsec}
In this section, we present a variant of the construction of Section
\ref{multconstsec}. This construction and the idea of the proof of the
following lemma has been taken from \cite[Prop. 5.4]{Osserman1}. We
fix integers $e_1>p$ and $1<e_2, e_3, e_4<p$ with $\gcd(e_1, p)=1$
such that $e_1+e_2+e_3+e_4$ is even and $e_3+e_4\leq
d:=(e_1+e_2+e_3+e_4-2)/2$.

\begin{lem}\label{addconstlem}
Let $k$ be an algebraically closed field of positive characteristic $p$.  
\begin{itemize}
\item[(a)] Let $f:\PP^1_k\to \PP^1_k$ be a cover with
ramification type $\tilde{\boldsymbol{C}}=(d; e_1, e_2, e_3-e_4)$. We
assume that the ramification points are $x=\infty, 0, 1, \rho$ and that $f(\infty)=\infty, f(0)=0$ and
$f(1)=f(\rho)=1$. Then for every $c\in \PP^1_k\setminus\{0, -1,
-\rho^{-p}, \infty\}$ the rational function
\[
f_c=f(x)+cx^p
\]
defines a cover of ramification type $\boldsymbol{C}=(d; e_1, e_2,
e_3, e_4)$.
\item[(b)] Conversely, assume that $g$ is a cover of ramification type 
  $\boldsymbol{C}=(d; e_1, e_2, e_3, e_4)$. Then there exists a $c\in
  k$ such that $g+cx^p$ has ramification type
  $\tilde{\boldsymbol{C}}=(d; e_1, e_2, e_3\text{-}e_4)$.
\end{itemize}
\end{lem}

\proof (a)  Let $f$ and $f_c$ be as in the statement of the lemma. We may
write
\[
f(x)=\frac{x^{e_2}f_1}{f_2},
\]
where $\deg(f_1)=d-e_2$ and $\deg(f_2)=d-e_1$. Moreover, the
polynomials $f_1$ and $f_2$  have simple zeros and are relatively prime. 
This implies that
\[
f_c(x)=\frac{x^{e_2}(f_1+cx^{p-e_2}f_2)}{f_2}.
\]
Therefore, the ramification index of $f_c$ in $x=0$ is $e_2$. Since
$e_1>p$ it follows also that the ramification index of $f_c$ in
$x=\infty$ is $e_1$. Similarly, the ramification
indices of $f_c$ in $x=1, \rho$ are $e_2, e_3$, respectively.

The equality
\[
\frac{\partial f_c}{\partial x}=\frac{\partial f}{\partial x}
\]
implies that $f_c$ is unramified outside $x=0, 1, \rho,\infty$. The
assumption on $c$ implies that the image of $x=0, 1, \rho,\infty$
under $f_c$ are all distinct and the statement of the lemma follows.

(b) Let $g$ be as in the statement of the lemma. Define $g_c=g+cx^p.$
Since $\partial g_c/\partial x=\partial g/\partial x\neq 0$ it follows
that $g_c$ is separable. Moreover, for all $c$ such that the image
under $g_c$ of the ramification points are pairwise distinct, the
ramification type of $g_c$ is still $\boldsymbol{C}=(d; e_1, e_2, e_3,
e_4)$. Assume that two of the ramification points, for example $x_3$
and $x_4$, have the same image under $g_c$. 
Then the ramification type is $\tilde{\boldsymbol{C}}=(d; e_1,
e_2, e_3\text{-}e_4)$.  The connectedness of the Hurwitz curve
$\HH_p(\boldsymbol{C})$ (Proposition \ref{4ptprop}.(b)) implies that
there exists a $c$ such that $g_c(x_3)=g_c(x_4)$, which proves (b).
\Endproof

The following proposition is a direct consequence of Lemma \ref{addconstlem}.

\begin{prop}\label{addconstprop} The assumptions being as above, assume additionally that $e_3\neq
e_4$. 
\begin{itemize}
\item[(a)]
We then have the equality
\[
h_p(d, e_1, e_2, e_3, e_4)=h_p(d; e_1, e_2, e_3\text{-}e_4).
\]
\item[(b)] If $h_p(d; e_1, e_2, e_3\text{-}e_4)>0$ then the Hurwitz
  curve $\HH_p(d, e_1, e_2, e_3, e_4)$ contains $h_p(d; e_1, e_2,
  e_3\text{-}e_4)$ irreducible components of genus $0$. Moreover, the restriction of
  the natural map $\pi: \HH_p(d, e_1, e_2, e_3, e_4)\to \PP^1_\lambda$
  which sends $[f]$ to its fourth branch point has degree $1$ on each
  of these components.
\end{itemize}
\end{prop}

\proof To prove (a), it is sufficient to show that nonisomorphic covers $f_i$ of
type $(d, e_1, e_2, e_3, e_4)$ give rise to nonisomorphic covers
under the construction of Lemma \ref{addconstlem}.

Let $f^i:\PP^1_k\to \PP^1_k$ be two nonisomorphic covers of type $(d,
e_1, e_2, e_3$-$e_4)$, and assume they are normalized as in the
statement of Lemma \ref{addconstlem}. The branch points of $f^i_c$ are
$\infty, 0, 1+c, 1+c\rho^p$. Normalizing the third branch point to $1$
yields the normalized cover $g^i_c(x):=f^i_c(x)/(1+c)$ with branch
points $\infty, 0, 1, (1+c\rho^p)/(1+c)=:\lambda_i$. The assertion
that the $g^i_c$ are nonisomorphic follows immediately from the
assumption that $e_3\neq e_4$. 

Statement  (b) follows immediately from the explicit expression
for the cover $f_c$ given in the proof of Lemma \ref{addconstlem}.
\Endproof

In the rest of this section, we discuss a concrete application of this
result to Hurwitz curves in positive characteristic.

\begin{lem}\label{3ptaddlem}
Let $p>3$ be a prime, and choose $e_1=p+2, e_2=3, 2\leq e_3< e_4<p$
with $e_3+e_4=p+1$. Put $d=(e_1+e_2+e_3+e_4-2)/2=p+2$.
We then have the inequality 
\[
h_p(d; e_1, e_2, e_3\text{-}e_4)\geq 1.
\] 
\end{lem}

\proof Let $f:\PP^1_k\to \PP^1_k$ be a cover of type $(d; e_1, e_2,
e_3\text{-}e_4)$, and assume that $f$ is normalized as in Lemma
\ref{addconstlem}. Then 
\begin{equation}\label{3ptaddeq}
f-1=c(x-1)^{e_3}(x-\rho)^{e_4}(x-a),
\end{equation}
for some $a\in \PP^1_k\setminus\{0, 1, \rho, \infty\}$. 

 Let $f$ be as in (\ref{3ptaddeq}). We want to determine which covers
 of this form define a cover of ramification type
 $\tilde{\boldsymbol{C}}=(d; e_1, e_2, e_3\text{-}e_4)$. 

The cover $f$ being ramified of order $3$ at $x=0$, we obtain the relations
\[
e_3a^2+2e_3a+2-e_3=0, \qquad \rho=\frac{3a+2-e}{e+1}.
\]
The coefficient $c$ is uniquely determined by $a$ and $\rho$ and the
condition $f(0)=0$.  A polynomial $f$ satisfying these conditions defines a cover of
ramification type $(d; e_1, e_2, e_3\text{-}e_4)$ if and only if the
values $0, 1, \infty, a$ and $ \rho$ are pairwise distinct. Our
assumptions on the $e_i$ imply in particular that $e_3\not\equiv
-1\pmod{p}$ and $d\equiv 2\pmod{p}$.  It then easily follows that $0, 1, \infty,
a,$ and $\rho$ are pairwise distinct if and only if
\begin{equation}\label{aexcepteq}
a\not\in \{(e_3-2)/2, (2e_3-1)/3, -1\}.
\end{equation}

If $a$ satisfies (\ref{aexcepteq}), then the ramification points $0,
1, \infty, a$ and $\rho$ are pairwise distinct for all $a$ satisfying
$e_3a^2+2e_3a+2-e_3=0$. Hence in this case we have $h_p(d; e_1, e_2,
e_3$-$e_4)=2$. 

Note that the integers $(e_3-2)/2, (2e_3-1)/3,$ and $-1$ are pairwise distinct, since
$e_3\not\equiv -1\pmod{p}$ by assumption. Therefore if $a\in
\{(e_3-2)/2, (2e_3-1)/3, -1\}$ the equation $e_3a^2+2e_3a+2-e_3=0$ has
exactly one solution $a$ for which the ramification points $0, 1,
\infty, a$ and $\rho$ are pairwise distinct. Therefore, in this case we find
$h_p(d; e_1, e_2, e_3$-$e_4)=1$.
\Endproof

The following result immediately follows from Proposition
\ref{addconstprop} and Lemma \ref{3ptaddlem}.

\begin{cor}\label{addconstcor}
Let 
$p>3$ be a prime, and choose $e_1=p+2, e_2=3, 2\leq e_3< e_4<p$
with $e_3+e_4=p+1$. Setting $d=(e_1+e_2+e_3+e_4-2)/2=p+2$, we then have the inequality 
\[
h_p(d; e_1, e_2, e_3, e_4)\geq 1.
\] 
\end{cor}

\begin{rem}
 The result of Proposition \ref{addconstprop}.(a) may also be
  deduced from Lemma \ref{addconstlem} by using deformation of
  admissible covers (see \cite[\S\ 2]{BO} and the references
  therein). However, this argument does not yield the information on
  the Hurwitz curve from Proposition \ref{addconstprop}.(b).
\end{rem}

\end{document}